\documentclass[12pt, reqno]{amsart}
\usepackage{amssymb, amsmath, amsthm, hyperref, graphicx}

\newtheorem{theorem}{Theorem}
\newtheorem{conj}[theorem]{Conjecture}

\begin{document}
\title{Around Efimov's differential test for homeomorphism}
\author{Victor Alexandrov}
\address{Sobolev Institute of Mathematics, Koptyug ave., 4, 
Novosibirsk, 630090, Russia and Department of Physics, 
Novosibirsk State University, Pirogov str., 2, Novosibirsk, 
630090, Russia}
\email{alex@math.nsc.ru}
\address{{}\hfill{September 20, 2020}}

\begin{abstract}
In 1968, N.\,V.~Efimov proved the following remarkable theorem:
\par
\textit{Let $f:\mathbb{R}^2\to\mathbb{R}^2\in C^1$ be such that $\det f'(x)<0$ 
for all $x\in\mathbb{R}^2$ and let there exist
a function  $a(x)>0$ and constants $C_1\geqslant 0$, $C_2\geqslant 0$ 
such that the inequalities
$|1/a(x)-1/a(y)|\leqslant C_1 |x-y|+C_2$
and
$|\det f'(x)|\geqslant a(x)|\operatorname{curl}f(x)|+a^2(x)$
hold true for all $x, y\in\mathbb{R}^2$. 
Then $f(\mathbb{R}^2)$ is a convex domain and $f$ maps $\mathbb{R}^2$ onto
$f(\mathbb{R}^2)$ homeomorphically.}
\par
Here $\operatorname{curl}f(x)$ stands for the curl of $f$ at $x\in\mathbb{R}^2$.
\par
This article is an overview of analogues of this theorem, its generalizations
and applications in the theory of surfaces, theory of global inverse functions, 
as well as in the study of the Jacobian Conjecture and the global asymptotic stability of dynamical systems.
\end{abstract}
\maketitle

\section{Introduction}
\label{sec:1}

One of the greatest achievements of 19th-century mathematics was the creation of Lobachevsky geometry, 
or, as it is now customary to say, hyperbolic geometry.
At that time, the issue of internal consistency of such geometry was of fundamental importance.
The first step in this direction was made by Eugenio Beltrami in 1868 in his article 
``Experience in the interpretation of non-Euclidean geometry'', see \cite{Be68}.
He proved that \textit{special domains of the Lobachevsky plane can be mapped isometrically onto 
suitable domains on the pseudosphere}.
Here the pseudosphere is, as usual, the surface in Euclidean 3-space formed by the rotation of 
the tractrix near its asymptote; and a domain is called special if it belongs to one of the following three classes: 
(i) the domains enclosed between two intersecting straight lines and a circle orthogonal to both of them, 
(ii) the domains enclosed between two divergent straight lines, their common perpendicular and the equidistant orthogonal to them, 
and (iii) the domains enclosed between two parallel straight lines and a horocycle orthogonal to both of  them.
However, the entire Lobachevsky plane cannot be mapped isometrically onto any domain on the pseudosphere.
Therefore, the problem of finding a regular surface in Euclidean 3-space which is isometric to the entire Lobachevsky 
plane was of fundamental importance for 19th-century mathematics.
However, in 1901, David Hilbert gave a negative solution to this problem, see \cite{Hi01}.
Namely, he proved that \textit{in Euclidean 3-space there is no complete regular surface 
with Gauss curvature $K=\operatorname{const} < 0$}. 
As a result, the proof of the consistency of Lobachevsky geometry went a different way, 
while the problem of how essential is the requirement of constancy of Gauss curvature in the Hilbert's theorem 
remained open for a long time.
Only in 1964 Nikola\u{\i} V.~Efimov completely removed this restriction by proving that

\begin{theorem}[{\cite[n.~1]{Ef64}}]\label{th:Ef64}
No surface can be $C^2$-immersed in Euclidean 3-space so as to be complete in the induced Riemannian metric, 
with Gauss curvature $K\leqslant\operatorname{const} < 0$.
\end{theorem}

In her famous article \cite[p. 474]{KM72}, Tilla Klotz Milnor writes that Efimov's proof of Theorem \ref{th:Ef64}
``is most ingenious, and does not depend upon sophisticated or modern techniques.''
After \cite{Ef64}, Efimov published several more articles related to Theorem \ref{th:Ef64}.
Among other things, he established that the condition $K\leqslant\operatorname{const} < 0$
is not the only obstacle for the immersibility of a complete surface of negative curvature.
In \cite{Ef66}, he showed that a rather slow change of Gauss curvature is another obstacle.
In all those numerous articles, he used to a large extent one and the same method based on the study 
of the spherical image of a surface. 
At that study, an essential role belongs to statements that, under some conditions, locally homeomorphic mapping 
$f:\mathbb{R}^2\to\mathbb{R}^2$ is, in the reality, a global homeomorphism and $f(\mathbb{R}^2)$ is a convex 
domain in $\mathbb{R}^2$. 
For example, in 1968 Efimov formulated and proved the following Theorems~\ref{th:Ef68a} and~\ref{th:Ef68b}:

\begin{theorem}[{\cite[Theorem~II]{Ef68}}]\label{th:Ef68a}
Let $f:\mathbb{R}^2\to\mathbb{R}^2$ be a continuously differentiable mapping such that 
$\det f'(x)<0$ for all $x\in\mathbb{R}^2$ and let there exist a function $a(x)>0$ and 
constants $C_1\geqslant 0$, $C_2\geqslant 0$ such that the inequalities
\begin{equation}\label{eq:1}
\biggl|\frac{1}{a(x)}-\frac{1}{a(y)}\biggr|\leqslant C_1 |x-y|+C_2
\end{equation}
and
\begin{equation}\label{eq:2}
|\det f'(x)|\geqslant a(x)|\operatorname{curl}f(x)|+a^2(x)
\end{equation}
hold true for all $x, y\in\mathbb{R}^2$. 
Then $f(\mathbb{R}^2)$ is a convex domain and $f$ maps $\mathbb{R}^2$ onto
$f(\mathbb{R}^2)$ homeomorphically.
\end{theorem}

Here $\operatorname{curl}f(x)=\partial f_2/\partial x_1 (x)-\partial f_1/\partial x_2 (x)$ 
stands for the curl of $f=(f_1, f_2)$ at $x=(x_1,x_2)\in\mathbb{R}^2$.

Note that if $1/a(x)$ has exponential growth (and, thus, the condition (\ref{eq:1}) is 
roughly violated) then the conclusion of Theorem \ref{th:Ef68a} 
may not hold even for potential mappings (i.\,e., such that $\operatorname{curl}f(x)\equiv 0$).
For example, for $f(x_1,x_2)=(e^{x_1}\sin x_2, e^{x_1}\cos x_2)$ 
we have $\operatorname{curl}f(x)\equiv 0$, $\det f'(x_1,x_2) = - e^{2x_1}<0$, 
$f(\mathbb{R}^2)=\mathbb{R}^2\diagdown \{(0,0)\}$ is nonconvex and every point from 
$f(\mathbb{R}^2)$ has infinite number of preimages.

\begin{theorem}[{\cite[Theorem~3]{Ef68}}]\label{th:Ef68b}
Let $f:\mathbb{R}^2\to\mathbb{R}^2$ be a continuously differentiable mapping such that 
$\det f'(x)<0$ for all $x\in\mathbb{R}^2$ and let there exists $a=\operatorname{const}>0$ 
such that for all $x\in\mathbb{R}^2$ the inequality
\begin{equation}\label{eq:3}
|\det f'(x)|\geqslant a|\operatorname{curl}f(x)|+a^2
\end{equation}
holds true. 
Then $f(\mathbb{R}^2)$ is either an infinite strip between parallel straight lines, or a half-plane, or the plane,
and $f$ maps $\mathbb{R}^2$ onto $f(\mathbb{R}^2)$ homeomorphically.
\end{theorem}

Note that the condition (\ref{eq:3}) of Theorem \ref{th:Ef68b} is satisfied provided
$|\det f'(x)|\geqslant \operatorname{const}> 0$ and 
$|\operatorname{curl}f(x)|\leqslant \operatorname{const}<+\infty$
hold true for all $x\in\mathbb{R}^2$.

This article is about the trace left by Theorems \ref{th:Ef68a} and \ref{th:Ef68b} in mathematics.
More precisely, we give an overview of the analogues of these theorems, their generalizations and 
applications over the past 40 years.
Sections \ref{sec:2}--\ref{sec:5} are devoted to presentation of results motivated by the theory of surfaces,
the theory of global inverse function, the Jacobian Conjecture, and the global asymptotic stability of dynamical systems,
respectively.
Of course, some results can be attributed to several sections; nevertheless, we present each result only once.
We also note that  
literature on the issues discussed in this article is extremely extensive
and we are not able to cover all the results known today.
The choice is limited by our personal preferences.

This article is an extended and updated version of the lecture  \cite{Al11} given by the author at 
the conference ``Metric geometry of surfaces and polyhedra'' dedicated to 
the centennial birthday anniversary of N.V.\,Efimov, held at Lomonosov Moscow State University in August 2010. 

Nikola\u{\i} Vladimirovich Efimov (May 31, 1910 -- August 14, 1982) was a prominent Soviet mathematician.
The reader can find information about his life and scientific heritage in his obituary \cite{AN83} and 
the book of memoirs about him \cite{TS14}.

\section{Results motivated by the theory of surfaces}\label{sec:2}

For completeness, we start with a mini-survey of results related to Theorem~\ref{th:Ef64}. 
                                       
A fully detailed exposition of Efimov's proof of Theorem \ref{th:Ef64} is given 
in the well-written article \cite{KM72}, which contributed a lot to popularity of Theorem \ref{th:Ef64} 
among western mathematicians.
Formally speaking, Theorem \ref{th:Ef64} is related to the general problem of immersing a Riemannian metric into a 
Euclidean space of suitable dimension.
Numerous refinements of this general problem are known;
extremely diverse technologies are used to solve them, see \cite{Gr17}, \cite{HH06}.
In reality, Theorem \ref{th:Ef64} has had the greatest impact on the field of mathematics known as geometry ``in the large.''
In this field, the emphasis is on non-local properties of isometric immersions (and embeddings) of manifolds into Euclidean 
spaces or other spaces under various restrictions on the curvature of the immersed manifold, and on connections between 
geometry and topology of immersed manifolds. 
The reader interested in the impact of Theorem \ref{th:Ef64} on various problems of geometry ``in the large'' is referred to
survey articles \cite{Am84}, \cite{Am19}, \cite{BS92}, \cite{Ro92}, \cite{RS11}. 
Below, we briefly mention only articles that are not included in these surveys.

In \cite{SX87},  the authors study the $n$-dimensional Efimov Conjecture: \textit{if $S$ is a complete isometrically immersed 
hypersurface in $\mathbb{R}^{n+1}$ whose Ricci curvatures are all negative, then these curvatures are not bounded above 
by a negative constant.} When $n=3$ the authors prove that this conjecture is correct in a very strong sense: 
the second fundamental form $A$ of such a hypersurface must satisfy $\inf |A|=0$. In the case $n>3$, they prove a partial result: 
if a counterexample to the conjecture exists, then its sectional curvatures must take every real value. 

In \cite{Pe89}, the author constructs an example of a complete analytic saddle surface in $\mathbb{R}^4$ with 
Gauss curvature $K\leqslant\operatorname{const} < 0$. 
This example shows that Theorem \ref{th:Ef64} does not generalize to saddle surfaces in $\mathbb{R}^n$, $n\geqslant 4$.
Recall that a complete surface in Euclidean space is called \textit{saddle} if it has no locally strictly
supporting hyperplane at any point; saddle surfaces have an intrinsic metric of nonpositive curvature.

In \cite{Sc95}, a theorem analogous to Theorem \ref{th:Ef64} is proved in the case of an ambient space of variable 
curvature with some extra conditions on its sectional curvature. 

In \cite{Sc01}, the author proves the following analogue of Theorem \ref{th:Ef64}: 
\textit{Let $S$ be a complete Riemannian surface with curvature $K\leqslant\operatorname{const} < 0$ 
such that $||\nabla K||\cdot |K|^{-3/2}$ is bounded. Then $S$ does not admit any isometric immersion neither into 
hyperbolic space $\mathbb{H}^3$ nor into spherical space $\mathbb{S}^3$.} Similar result is also proved for Lorentz 3-space 
$\mathbb{H}^3_1$ of constant curvature $-1$.

In \cite{GMT15a}, the authors study complete surfaces in $\mathbb{R}^3$ with Gauss curvature satisfying the inequality  
$K\leqslant\operatorname{const} < 0$ in a neighborhood of infinity. They prove that such a surface is topologically a finitely 
punctured compact surface, its surface area is finite, and each puncture is a cusp extending to infinity, asymptotic to a ray.
Similar results on complete surfaces in non-Euclidean space forms are obtained in  \cite{GMT15b}.

In \cite{Me18}, the author studies non-compact connected surfaces $S$ with compact boundary and proves that
there is no complete isometric immersion of $S$ into $\mathbb{R}^3$ satisfying that $\int_S |K | = +\infty$ and 
$K\leqslant\operatorname{const} < 0$. In particular, he shows that Theorem \ref{th:Ef64} holds true for complete Hadamard 
immersed surfaces (i.\,e., immersed complete  simply  connected  Riemannian manifolds with non positive sectional curvature), 
whose Gauss curvature $K$ is bounded away from zero outside a compact set.

In \cite{Ch00}, \cite{Ch03}, and \cite{Ch06}, it is shown that Theorem \ref{th:Ef64} allows us to better understand 
the solutions of the Einstein--Maxwell equation of gravitation and electromagnetism in source-free space. 

Summing up this mini-survey of results related to Theorem \ref{th:Ef64}, we can say that it has had and continues to 
have a significant impact on geometry ``in the large,''
some problems closely related to it still remain to be answered; also they attract the attention of young researchers. 
The latter is evidenced, for example, by the Master's Dissertation \cite{Ar19}. 
The most well-known still open problem closely related to Theorem \ref{th:Ef64} reads as follows:

\begin{conj}[John Milnor, \cite{KM72}]\label{conj:Milnor}
Suppose $S$ is a complete, umbilic free
surface, $C^2$-immersed in $\mathbb{R}^3$ so that the sum of the squares of the principal
curvatures on $S$ is bounded away from zero. Then either $K$ changes sign on
$S$, or else $K=0$.
\end{conj}
Some results related to Conjecture \ref{conj:Milnor} may be found in \cite{To95}, \cite{To96}.

Now we turn to a discussion of results which are related to Theorems \ref{th:Ef68a} and \ref{th:Ef68b},
were obtained in connection with the theory of surfaces, and were not presented in the 
survey articles \cite{Am84}, \cite{Am19}, \cite{BS92}, \cite{Ro92}, \cite{RS11}. 

Under the conditions of Theorem \ref{th:Ef68b}, $f(\mathbb{R}^2)$ is either an infinite strip between 
parallel straight lines, or a half-plane, or the plane. 
Results clarifying which of these possibilities may actually occur can be found in 
the articles \cite{Ka70} and \cite {Ge70}, book \cite[\S\,31]{BVK73}, and Ph.~D. thesis \cite{Pe90}.
Let us describe them in some detail.

The case, when $f(\mathbb{R}^2)$ is the plane, is, obviously, realized by the linear mapping $f_1(x_1,x_2)=x_2$, $f_2(x_1,x_2)=x_1$.

The case, when $f(\mathbb{R}^2)$ is a half-plane, is realized by the mapping $f:\mathbb R^2\to\mathbb R^2$ defined 
by the formulas found in \cite{Ka70}:
\begin{equation}\label{eq:4}
f_1(x_1,x_2) = \log\bigl(x_1+\sqrt{x_1^2+ e^{-2x_2}}\bigr) + x_2, \quad
f_2(x_1,x_2) = \sqrt{x_1^2 + e^{-2x_2}}.
\end{equation}
In fact, it is easy to check directly that for this mapping we have 
\begin{equation}\label{eq:5}
\det f'(x)= \operatorname{const}<0, \quad \operatorname{curl}f(x)= 0,
\end{equation} 
and $f(\mathbb R^2)$ is the half-plane $f_2>0$. 
The formulas (\ref{eq:4}) look complicated, but they can be obtained in the following natural way proposed in \cite{Ka70}.
Note that the conditions (\ref{eq:5}) will certainly be satisfied
if $f=(f_1,f_2)$ is potential (i.\,e., if there exists a function
$\varphi:\mathbb{R}^2\to\mathbb{R}$, called a potential, such that  
$f_j=\partial\varphi/\partial x_j$ for $j=1,2$) and the potential $\varphi$
satisfies the Monge--Amp\`{e}re equation
\begin{equation*}
\frac{\partial^2 \varphi}{\partial x_1^2}\frac{\partial^2 \varphi}{\partial x_2^2}-
\biggl(\frac{\partial^2 \varphi}{\partial x_1\partial x_2}\biggr)^2=-1.
\end{equation*}
Then using the classical parametric representation for solutions of the Monge--Amp{\`e}re equation 
\cite[Section 476, example 4]{Go64}, we obtain a large class of non-trivial examples of mappings $f:\mathbb{R}^2\to\mathbb{R}^2$ 
satisfying the conditions (\ref{eq:5}) for which $f(\mathbb{R}^2)$ is a half-plane.
The formulas (\ref{eq:4}) define just one mapping from that class.

The case, when $f(\mathbb{R}^2)$ is an infinite strip between parallel straight lines
was studied in \cite{Ka70} and \cite{Ge70}. 
More precisely, in \cite{Ka70} it is proved that this case is impossible provided that the conditions (\ref{eq:5}) are fulfilled.
Moreover, in \cite{Ka70} it is conjectured that this case can never occur under the conditions of Theorem \ref{th:Ef68b}.
This conjecture remains open to this day. Partial results on this conjecture are obtained 
in \cite{Ge70}, where it is shown that $f(\mathbb{R}^2)$ cannot be a strip
provided that $\operatorname{curl}f(x)= 0$ and $\det f'(x)=-g^2(x)$ with some function $g$,
which satisfies the inequality $g(x)\geqslant\operatorname{const}>0$ and is either a convex function or a
polynomial.

We proceed by describing the results obtained by G.\,Ya.~Perel'man in his Ph.~D. dissertation \cite{Pe90}.
The dissertation is devoted to the study of saddle surfaces in $\mathbb{R}^n$, $n\geqslant 3$, i.\,e., 
to generalization and strengthening 
of theorems proved by Efimov in \cite{Ef64}, \cite{Ef66}, \cite{Ef68} and, in particular, of Theorem \ref{th:Ef64}.
Among other results, G.\,Ya.~Perel'man proved the following theorem which, as far as we know, was never published in a
peer reviewed journal: 

\begin{theorem}[{\cite[p.~78]{Pe90}}]\label{th:Pe90}
Let $f:\mathbb{R}^2\to\mathbb{R}^2$ be a potential orientation-reversing global diffeomorphism of class $C^1$. 
Then $f(\mathbb{R}^2)$ is a domain from the following list:
$(1)$ a triangle;
$(2)$ a quadrangle;
$(3)$ the plane;
$(4)$ a half-plane;
$(5)$ a strip between two parallel lines;
$(6)$ an angle;
$(7)$ a half-strip;
$(8)$ an angle with a cut off vertex;
$(9)$ a half-strip with a cut off vertex.
\end{theorem}
Note that the conditions of Theorem~\ref{th:Pe90} do not control the decay of $\det f'(x)$ at infinity at all,
i.\,e., no condition similar to the inequalities (\ref{eq:2}) or (\ref{eq:3}) is imposed.
The proof of Theorem~\ref{th:Pe90} is a clever modification of the Efimov's method, proposed in 
\cite{Ef64}, \cite{Ef66}, \cite{Ef68}.

We complete this Section with a description of an extension of Theorem~\ref{th:Ef64} to surfaces in the 
three-dimensional Heisenberg group $Nil$ published in \cite{BP11}.
Recall that the Heisenberg group $Nil$ is $\mathbb{R}^3$ endowed with the metric
\begin{equation*}
ds^2=dx^2+dy^2+\bigl(\tfrac12 ydx-\tfrac12 xdy +dz
\bigr)^2.
\end{equation*} 
Here $x,y,z$ are global Cartesian coordinates in $\mathbb{R}^3$. 
It is known that, in the Heisenberg group, surfaces with Gauss curvature bounded above by a negative constant do exist.
However, the following analogue of Theorem~\ref{th:Ef64} holds true for surfaces, which are explicitly given in 
coordinates $(x, y)$:

\begin{theorem}[{\cite[Theorem 3]{BP11}}]\label{th:BP}
In $Nil$, on any $C^2$-surface $z=f(x,y)$ explicitly given in the entire plane $(x,y)$,
the greatest lower bound of the absolute value of the Jacobian of Gauss mapping is zero.
\end{theorem}
In Theorem~\ref{th:BP}, the replacement of Gauss curvature by the Jacobian of Gauss mapping looks quiet natural, and
Gauss mapping of a surface in $Nil$ is defined as Gauss mapping of a surface in a Lie group.
For our exposition it is also important that Theorem~\ref{th:Ef68b} is used essentially in the proof of Theorem~\ref{th:BP}. 
Indeed, the key point of the proof of Theorem~\ref{th:BP} is that the mapping $F:\mathbb{R}^2\to\mathbb{R}^2$ defined as 
the composition of Gauss mapping and the projection of the sphere from the origin onto its tangent plane 
at the point $(0,0,1)$ satisfies the conditions of Theorem \ref{th:Ef68b} and, thus, is a homeomorphism
of $\mathbb{R}^2$ onto a convex domain $F(\mathbb{R}^2)$ in~$\mathbb{R}^2$.

\section{Results motivated by theorems on global inverse function}\label{sec:3}

Theorems on global inverse function play an important role in many fields of Analysis and its applications,
including mathematical economics.
The interested reader may find a broad range of such theorems, e.\,g., in \cite {Pl74}, \cite{Pa83}, \cite{Gr17}, and \cite{JLM19}.
Geometric applications may be found, e.\,g., in the book \cite{Sa08} devoted to the study of 
isometric immersions and embeddings of locally Euclidean metrics using various theorems on the univalence of mappings, and 
in the survey paper \cite{Bo10} devoted to the study of univalent solutions of the uniformly elliptic Beltrami
equation realizing a homeomorphic  quasiconformal mapping of the complex plane with the assigned measurable complex dilatation.

To the best of our knowledge, the first result on existence of a global inverse function was obtained 
by Jacques Hadamard in \cite{Ha06} for mappings $\mathbb{R}^2\to\mathbb{R}^2$. 
Nowadays the following statement of the global inverse function theorem is considered standard:
\begin{theorem}[{\cite[Theorem 3.2]{Pl74}}]\label{th:Pl74}
Let $X$ and $Y$ be Banach spaces, $f: X \to Y\in C^1(X)$ and 
$f'(x)$ is invertible for all $x\in X$. If
\begin{equation}\label{eq:6}
\int\limits_{0}^{+\infty} \inf\limits_{\|x\|\leqslant t}{\|f'(x)^{-1}\|^{-1}}\,dt=+\infty,
\end{equation}
then $f$ is a diffeomorphism of $X$ onto $Y$.
\end{theorem}
Note that the condition (\ref{eq:6}) of Theorem~\ref{th:Pl74} controls the growth of the norm $\|f'(x)^{-1}\|$
at infinity and is obviously fulfilled in the case, when
\begin{equation}\label{eq:7}
\|f'(x)^{-1}\|\leqslant A=\operatorname{const}<+\infty
\end{equation}
for all $x\in X$. 

Hadamard's theorem and its modern versions, e.\,g., Theorem \ref{th:Pl74}, admit numerous generalizations 
to nonsmooth mappings. As an example, we mension the following theorem
\begin{theorem}[{\cite[Theorem 2.2]{AZ19}}]\label{th:AZ19}
Let $n\geqslant 1$, $\alpha >0$, and let $f:\mathbb{R}^n\to\mathbb{R}^n$ be continuous, locally injective 
and $\alpha$-covering at each $x\in\mathbb{R}^n$ 
$($the latter means that, for every $\varepsilon>0$, there exists $r\in (0,\varepsilon ]$ such that 
$B(f(x), \alpha r)\subset f(B(x, r))$, where $B(x,r)$ stands for the ball of radius $r$ centered at $x)$. 
Then $f$ is a global homeomorphism and $f^{-1}$ is Lipschitz continuous with the constant $\alpha^{-1}$.
\end{theorem}

In \cite{Al90}, it is proved that the following theorem is equivalent to Theorem~\ref{th:Ef68a}:
\begin{theorem}[{\cite[Theorem 4]{Al90}}]\label{th:Al90}
Let $f:\mathbb{R}^2\to\mathbb{R}^2$ be a continuously differentiable mapping such that 
$\det f'(x)<0$ for all $x\in\mathbb{R}^2$; and let there exist a function $a(x)>0$ and 
constants $C_1\geqslant 0$, $C_2\geqslant 0$ such that the inequality $(\ref{eq:1})$
holds true for all $x, y\in\mathbb{R}^2$
and, for all $x\in\mathbb{R}^2$, the inequality
\begin{equation}\label{eq:8}
|\mu_2(x)|\geqslant |\mu_1(x)|\geqslant a(x)  
\end{equation}
holds true. Here $\mu_1(x)$ and $\mu_2(x)$ are the eigenvalues of $f'(x)$.
Then $f$ maps $\mathbb{R}^2$ onto $f(\mathbb{R}^2)$ homeomorphically, and $f(\mathbb{R}^2)$ is convex.
\end{theorem}
Note that the inequalities (\ref{eq:8}) control the decrease of the spectral radius $\rho\bigl(f'(x)\bigr)$ of 
$f'(x)$ at infinity and are obviously fulfilled in the case, when
\begin{equation*}
|\mu_2(x)|\geqslant |\mu_1(x)|\geqslant a=\operatorname{const}>0
\end{equation*}
for all $x\in \mathbb{R}^2$. As usual, the spectral radius of a linear mapping is equal to the maximum of the modules 
of its eigenvalues. 
In other words, we can say that the inequalities (\ref{eq:8}) control the growth of the spectral radius 
$\rho\bigl(f'(x)^{-1}\bigr)$ of $f'(x)^{-1}$
at infinity and are obviously fulfilled in the case, when
\begin{equation}\label{eq:9}
\rho\bigl(f'(x)^{-1}\bigr)\leqslant A=1/a=\operatorname{const}<+\infty
\end{equation}
for all $x\in \mathbb{R}^2$. Comparing (\ref{eq:7}) and (\ref{eq:9}), we arrive at the following conjecture:

\begin{conj}[{\cite[p. 199]{Al90}, \cite[p. 104]{Al91}}]\label{conj:Alexandrov}
Let $n\geqslant 2$, $f:\mathbb{R}^n\to\mathbb{R}^n$ be a continuously differentiable mapping such that 
$\det f'(x)\neq 0$ and  
$$\rho\bigl(f'(x)^{-1}\bigr)\leqslant A$$ 
for all $x\in\mathbb{R}^n$ with some constant $A<+\infty$.
Then $f(\mathbb{R}^n)$ is a convex domain in $\mathbb{R}^n$ and $f$ is injective.
\end{conj}

Note that Theorem~\ref{th:Al90} gives the affirmative answer to Conjecture~\ref{conj:Alexandrov} in the case, when
$n=2$ and $\det f'(x)<0$.
Note also that Conjecture~\ref{conj:Alexandrov} is in agreement with Theorem~\ref{th:Pl74}
in the following sense. Let $n\geqslant 2$, $f:\mathbb{R}^n\to\mathbb{R}^n$ 
be a continuously differentiable mapping,
$f'(x)$ be invertible and $\|f'(x)^{-1}\|\leqslant A=\operatorname{const}<+\infty$ for all $x\in \mathbb{R}^n$.
Then, according to Theorem~\ref{th:Pl74}, $f$ is a diffeomorphism of $\mathbb{R}^n$ onto $\mathbb{R}^n$.
At the same time, we see that such $f$ satisfies both the conditions of Conjecture~\ref{conj:Alexandrov}
(since $\rho\bigl(f'(x)^{-1}\bigr)\leqslant \|f'(x)^{-1}\|\leqslant A$) and its conclusion 
(since $f$ is injective and $f(\mathbb{R}^n)=\mathbb{R}^n$ is convex).

Conjecture~\ref{conj:Alexandrov} is apparently very natural and was formulated by different researchers in connection 
with problems from various fields of mathematics.
For example, in the article \cite{Al90} and conference proceedings \cite{Al91} it
was formulated in the form given above and was motivated by its connection with Efimov's Theorem \ref{th:Ef68a}. 
In \cite[Conjecture 2.1]{CM98}, it was formulated without mentioning that  $f(\mathbb{R}^n)$ is convex
and was motivated by its connection with the Jacobian Conjecture of Ott-Heinrich Keller \cite{Ke39};
see Section~\ref{sec:4} below for more details.
Only during the preparation of this survey article, we found that in fact for the first time 
Conjecture~\ref{conj:Alexandrov} was formulated in \cite{Ko77} in connection with Efimov's Theorem~\ref{th:Ef68a}.
Moreover, the following theorem is proved there:

\begin{theorem}[{\cite[p. 224]{Ko77}}]\label{th:Ko77}
Let $f:\mathbb{R}^2\to\mathbb{R}^2$ be a continuously differentiable mapping such that 
the eigenvalues $\mu_1(x)$, $\mu_2(x)$ of $f'(x)$ satisfy everywhere the conditions:

$(1)$ $\mu_1(x)$ and $\mu_2(x)$ are real;

$(2)$ $\mu_1(x)\neq\mu_2(x)$;

$(3)$ $0<a\leqslant \min\{|\mu_1(x)|, |\mu_2(x)|\}$ with some $a=\operatorname{const}$.

\noindent{Then} $f$ is injective and  $f(\mathbb{R}^2)$ is convex.
\end{theorem}
Note that Theorem~\ref{th:Ko77} is applicable to both orientation-reversing and orienta\-tion-preserving mappings $f$.
The proof of Theorem~\ref{th:Ko77} follows closely Efimov's proof of his Theorem~\ref{th:Ef68a}.

To conclude this Section, we point out another natural problem related to Conjecture~\ref{conj:Alexandrov}, 
which remains open so far. 
It reads as follows: \textit{find an estimation of the radius of the ball, where 
the inverse mapping $f^{-1}$ exists, in terms of the spectral radius $\rho\bigl(f'(x)^{-1}\bigr)$.} 
The estimation is expected to be similar to that in the following theorem:

\begin{theorem}[{\cite[Theorem IIA]{Jo68}}]\label{th:Jo68}
Let $X$ and $Y$ be Banach spaces, 
$B$ be the ball of center $x_0$ and radius $r$ in $X$,
and $f:B\to Y$ be a continuously differentiable mapping, whose derivative $f'(x)$ is invertible for all $x\in B$.
Then the inverse mapping $f^{-1}$ exists in the ball of center $f(x_0)$ and radius
\begin{equation*}
R=\int\limits_{0}^{r} \inf\limits_{\|x-x_0\|\leqslant t}{\|f'(x)^{-1}\|^{-1}}\,dt.
\end{equation*}
\end{theorem}
Theorem~\ref{th:Jo68} generalizes Theorem~\ref{th:Pl74} and is quoted here in our notation.

\section{Results motivated by the Jacobian Conjecture}\label{sec:4}

A mapping $f:\mathbb{C}^n\to\mathbb{C}^n$ is said to be polynomial if $f=(f_1, \dots, f_n)$ and each 
component $f_i$ of $f$ is a polynomial. 
A polynomial mapping $f$ is said to be a polynomial automorphism of $\mathbb{C}^n$ if $f$ is one-to-one,
$f(\mathbb{C}^n)=\mathbb{C}^n$, and  $f^{-1}:\mathbb{C}^n\to\mathbb{C}^n$ is also a polynomial mapping. 

\begin{conj}[Jacobian Conjecture, \cite{Ke39}]\label{conj:jacobian}
For every positive integer $n$ and every polynomial mapping $f:\mathbb{C}^n\to\mathbb{C}^n$, whose Jacobian 
is a non-zero constant, $f$ is a polynomial automorphism of $\mathbb{C}^n$.
\end{conj}
The Jacobian Conjecture was first formulated by Ott-Heinrich Keller in 1939 but is still unproved, even for $n=2$. 
Its history, many references, and some partial results can be found in \cite{Es00}.
Stephen Smale included the Jacobian Conjecture in his list of eighteen problems for the 21st century 
\cite[Problem 16]{Sm00}.
Note that various additional conditions on $f$ are known under which the Jacobian Conjecture is certainly valid.
One of such additional conditions is the requirement that $f$ is injective.
An elementary proof of this fact is given in~\cite {Ru95}.
This explains the interest of specialists on the Jacobian Conjecture in conditions guaranteeing the injectivity 
of mappings. We indicate one of their conjectures in our notation:

\begin{conj}[M.~Chamberland, {\cite[Conjecture 2.1]{CM98}}]\label{conj:Chamberland}
Let $n\geqslant 2$, $f:\mathbb{R}^n\to\mathbb{R}^n$ be a continuously differentiable mapping such that 
$\det f'(x)\neq 0$ and  
\begin{equation}\label{eq:10}
\rho\bigl(f'(x)^{-1}\bigr)\leqslant A
\end{equation} 
for all $x\in\mathbb{R}^n$ with some
constant $A<+\infty$,
where $\rho\bigl(f'(x)^{-1}\bigr)$ is the spectral radius of $f'(x)^{-1}$.
Then $f$ is injective.
\end{conj}
Note that the condition~(\ref{eq:10}) can be replaced (and is replaced in~\cite{CM98}) by the equivalent 
requirement that there exists an $\varepsilon>0$ such that $|\lambda|\geqslant \varepsilon$ for all 
the eigenvalues $\lambda$ of $f'(x)$ for all $x\in\mathbb{R}^n$.
Note also that, strictly speaking, Conjecture~\ref{conj:Chamberland} 
is weaker than Conjecture~\ref{conj:Alexandrov}, since it does not state that $f(\mathbb{R}^n)$ is convex.
However, Conjecture~\ref{conj:Chamberland} also remains open.
We consider the following theorem as closely related to Conjecture~\ref{conj:Chamberland}:

\begin{theorem}[{\cite[Theorem 1.1]{CM98}}]\label{th:CM98}
Let $n\geqslant 2$ and $f:\mathbb{R}^n\to\mathbb{R}^n$ be a continuously differentiable mapping.
Suppose there exists an $\varepsilon>0$ such that $|\mu|\geqslant \varepsilon$ for all the 
eigenvalues  $\mu$ of $f'(x)f'(x)^T$ for all $x\in\mathbb{R}^n$.
Then $f$ is injective.
\end{theorem}
To prove Theorem~\ref{th:CM98}, the authors use the so-called Mountain Pass Lemma due to \cite{AR73},
i.\,e., they use techniques which differs very much from those used by Efimov to
prove Theorems~\ref{th:Ef68a} and~\ref{th:Ef68b}.

\section{Results motivated by the global asymptotic stability of dynamical systems}\label{sec:5}

Conjectures similar to Conjectures~\ref{conj:Alexandrov} and~\ref{conj:Chamberland}
appeared also in the study of the global asymptotic stability of dynamical systems in the plane.
Namely, in this field of mathematics several authors conjectured that some restrictions 
on the set of all eigenvalues of $f'(x)$ (like separation from zero)
already imply the injectivity of a $C^1$-mapping $f:\mathbb{R}^n\to\mathbb{R}^n$.
In this Section we present some conjectures and theorems of this sort.

The problem of the global asymptotic stability of dynamical systems in the plane can be formulated as follows:
\textit{Consider an autonomous system in the plane
\begin{equation}\label{eq:11}
\begin{cases}
\dot x_1 &= f_1(x_1,x_2),\\
\dot x_2 &= f_2(x_1,x_2),
\end{cases}
\end{equation}
for which  $(x_1, x_2)=(0,0)$ is an equilibrium point $($i.\,e., $f_1(0,0)=f_2(0,0)=0)$ 
and the eigenvalues of the Jacobi matrix $f'(x)$ of the mapping $f=(f_1,f_2)$
have negative real parts at any point $x$ on the plane.
Is it then true that the trivial solution $(x_1, x_2)=(0,0)$ is
globally asymptotically stable, i.\,e., that any solution of the equations $(\ref{eq:11})$
tends to the point $(0,0)$ as $t\to\infty$?}
This problem plays an important role in the theory and applications of dynamical systems, see, e.\,g., \cite{Chi06}.
According to the famous Lyapunov theorem, an autonomous system of linear differential equations
$\dot{y}=Ay$, $y\in\mathbb R^n$, $n\geqslant 2$,
is globally asymptotically stable if all the eigenvalues of the constant matrix~$A$ have negative real parts.
Thus, the problem of the global asymptotic stability answers the question of whether it is 
possible to make a conclusion about the global asymptotic stability of any solution of the nonlinear system~(\ref{eq:11}) 
as soon as the global asymptotic stability of any solution of the linearized system 
$\dot{y}=f'(x)y$ is established for all $x\in\mathbb R^2$.

The problem of the global asymptotic stability of dynamical systems can be posed 
in $\mathbb R^n$ for arbitrary $n$. 
It is known as the Markus--Yamabe Conjecture and is known to be false for $n\geqslant 3$,
see, e.\,g., \cite{CEG97}.

In \cite{Ol63}, a close connection is discovered between the problem of the global asymptotic stability in the plane
and injectivity of some mappings; namely, the following theorem is proved there:

\begin{theorem}[{\cite{Ol63}}]\label{th:Ol63}
The following two statements are equivalent to each other: 
$(1)$ the problem of the global asymptotic stability in the plane has the affirmative answer{\rm;}
$(2)$ every continuously differentiable mapping $f:\mathbb{R}^2\to\mathbb{R}^2$, 
the eigenvalues of the Jacobi matrix $f'$ of which 
everywhere have negative real parts, is injective.
\end{theorem}

In 1994--1995, the affirmative answer to the problem of the global asymptotic stability
in the plane was obtained in \cite{Gl94}, \cite{Gl95}, \cite{Fe95}, and \cite{Gu95}.
The proofs proposed in these articles rely on Theorem~\ref{th:Ol63}.
However, interest in the problem of the global asymptotic stability in the plane does not disappear; 
on the contrary, more and more new proofs appear, see, e.\,g.,~\cite{CHQ01}.
Moreover, specialists in different fields of mathematics invent and study various differential 
tests for injectivity of mappings related to Theorems~\ref{th:Ef68a} and~\ref{th:Ef68b}, 
and Conjectures~\ref{conj:Alexandrov} and~\ref{conj:Chamberland}.
We consider the following theorem as a very important contribution to this field of research:

\begin{theorem}[{\cite[Theorem A]{CGL02}}]\label{th:CGL02}
Let $f:\mathbb{R}^2\to\mathbb{R}^2$ be  a continuously differentiable mapping, and let
$\operatorname{Spec}(f)$ be the set of all eigenvalues of the Jacobian matrix $f'(x)$ 
for all $x\in\mathbb{R}^2$. Suppose that
$\operatorname{Spec}(f)\cap(-\varepsilon, \varepsilon)= \varnothing$
for some $\varepsilon >0$.
Then $f$ is injective. 
\end{theorem}
The key point in the proof of Theorem~\ref{th:CGL02} is the notion of a half-Reeb component of a planar vector field.
In particular, the authors use techniques which is very different from those used by Efimov to
prove Theorems~\ref{th:Ef68a} and~\ref{th:Ef68b}.
Note that, for $n=2$, Theorem~\ref{th:CGL02} shows that a statement even stronger than 
Conjecture~\ref{conj:Chamberland} holds true. 
Namely, in this case for injectivity of $f$ it is sufficient to require that $\operatorname{Spec}(f)$ 
does not intersect with a sufficiently small interval $(-\varepsilon, \varepsilon)$ of the real axis $\mathbb{R}$ only, 
not with the entire disk of radius $\varepsilon$ centered at 0 in the complex plane $\mathbb{C}$.
In \cite{FGR04b}, it is shown that, in Theorem~\ref{th:CGL02}, the condition 
$\operatorname{Spec}(f)\cap(-\varepsilon, \varepsilon)= \varnothing$
can be replaced by the weaker condition 
$\operatorname{Spec}(f)\cap[0, \varepsilon)= \varnothing$.
On the other hand in~\cite{SX96}, it is proved that there exist integers $n>2$ and non-injective polynomial
mappings $f:\mathbb{R}^n\to\mathbb{R}^n$ with $\operatorname{Spec}(f)\cap[0, +\infty)= \varnothing$.
Nevertheless, some partial extensions of Theorem~\ref{th:CGL02} to $\mathbb{R}^n$, $n>2$, 
are obtained in~\cite{SX96} and~\cite{FGR04a}.

At last, we would like to mention the article~\cite{Ca00}, where the author explores whether the $C^1$-mappings 
with a unipotent (i.\,e., all eigenvalues are 1) Jacobian matrix are injective. 
He extends some results known in the polynomial case to the $C^1$ context, resolves some special cases, and
discusses connections between the problem of the injectivity of the $C^1$-mappings with a unipotent Jacobian matrix 
and the Jacobian Conjecture~\ref{conj:jacobian}, Chamberland's Conjecture~\ref{conj:Chamberland}, and 
a few additional ones that have arisen in connection with the study of the
Markus--Yamabe Conjecture and some its discrete analogue.
One of conjectures discussed in~\cite{Ca00} reads as follows:
\begin{conj}[$C^1$ Fixed Point Conjecture, {\cite[Conjecture 5]{Ca00}}]\label{conj:fixed-point}
If $f:\mathbb{R}^n\to\mathbb{R}^n$ is $C^1$ with $f(0)=0$, and the eigenvalues of $f'(x)$
have absolute value less than $1$ at every point $x\in\mathbb{R}^n$, then $0$ is the unique fixed point of $f$.
\end{conj}
According to~\cite{Ca00}, being applied to polynomial mappings 
Conjecture~\ref{conj:fixed-point} is equivalent to the Jacobian Conjecture~\ref{conj:jacobian}.
\section{Conclusion}\label{sec:6}
Theorems deducing the injectivity of mappings $f:\mathbb{R}^n\to\mathbb{R}^n$ from some properties 
of $f'(x)$ appear in different branches of mathematics.
Efimov's Theorems~\ref{th:Ef68a} and \ref{th:Ef68b} are among the first results of this type.
They directly stimulated the studies of some mathematicians, or later
turned out to be closely related to findings of other mathematicians not familiar with the works of Efimov.
Open problems related to Theorems~\ref{th:Ef68a} and \ref{th:Ef68b} continue to attract the attention of mathematicians.

\end{document}